\documentclass[12pt,a4paper]{article}
\usepackage{amsmath,amssymb,fancyhdr,url}
\usepackage{graphicx}
\usepackage{algorithmic}
\usepackage{algorithm}
\usepackage{listings}
\usepackage{graphicx}
\usepackage{amsfonts}
\usepackage{graphicx,color}
\usepackage{graphpap,color}
\usepackage{latexsym}
\usepackage{float}
\usepackage{setspace}
\usepackage{algorithmic}
\usepackage{algorithm}
\usepackage{listings}
\usepackage{rotating}
\usepackage{rotfloat}
\begin{document}
\title{A Switching Approach to Avoid Breakdown in Lanczos-type Algorithms}
\author{Muhammad Farooq and Abdellah Salhi}

\date{}
\maketitle
\bibliographystyle{plain}

\begin{abstract}
Lanczos-type algorithms are well known for their inherent instability. They typically breakdown when relevant orthogonal polynomials do not exist. Current approaches to avoiding breakdown rely on jumping over the non-existent polynomials to resume computation. This jumping strategy may have to be used many times during the solution process. We suggest an alternative to jumping which consists in switching between different algorithms that have been generated using different recurrence relations between orthogonal polynomials. This approach can be implemented as three different strategies: ST1, ST2, and ST3. We shall briefly recall how Lanczos-type algorithms are derived. Four of the most prominent such algorithms namely $A_4$, $A_{12}$, $A_5/B_{10}$ and $A_5/B_8$ will be presented and then deployed in the switching framework. In this paper, only strategy ST2 will be investigated. Numerical results will be presented.
\end{abstract}



\newpage
\bigskip
\section{Introduction}
Lanczos-type methods for solving SLE's are based on the theory of FOP's. All such methods are implemented via some recurrence relationships between polynomials $P_{k}(x)$ represented by $A_{i}$ or between two adjacent families of orthogonal polynomials $P_{k}(x)$ and $P^{(1)}_{k}(x)$ represented by $A_{i}$ and $B_{j}$ as described in \cite{94:Baheux,95:Baheux,10:Farooq}. The coefficients of the various recurrence relationships between orthogonal polynomials are given as ratios of scalar products. When a scalar product in a denominator vanishes, then a breakdown occurs in the algorithm and the process normally has to be stopped. Equivalently, the breakdown is due to the non-existence of some orthogonal polynomial or polynomials. So, an important issue is how to continue the solution process in such a situation and arrive at a useable result. Several procedures for that purpose appeared in the literature in the last few decades. It has been shown, for instance, that it is possible to jump over non-existing polynomials, \cite{94:Baheux,93:Brezinski}; breakdown-free algorithms were thus obtained. The first attempt in this regard was the look-ahead Lanczos algorithm, \cite{85:Parlett}. Other procedures for avoiding breakdown are also proposed in \cite{92:Gut,91:Brezinski,92:Sadok,94:Baheux,93:Brezinski,99:Brezinski,97:Graves,93:Gut,79:Parlett}. However, they all have their limitations
including the possibility of calling the procedure for remedying the breakdown, more than once. In the following, we suggest an alternative to jumping over missing polynomials by switching between different variants of the Lanczos algorithm.

\section{The Lanczos approach}
We consider a linear system of equations,
\begin{equation}
\textit{A}\textbf{x}=\textbf{b},
\end{equation}
\noindent where $\textit{A}\in \textit{R}^{n\times n}$,
$\textbf{b}\in \textit{R}^{n}$ and $\textbf{x}\in
\textit{R}^{n}$.

Let $\textbf{x}_{0}$ and $\textbf{y}$ be two arbitrary vectors in $\textit{R}^{n}$ such that $\textbf{y}\neq0$. The Lanczos method, \cite{52:Lanczos} consists in constructing a sequence of vectors $\textbf{x}_{k}\in \textit{R}^n$ defined as follows, \cite{02:Brezinski,95:Baheux}
\begin{equation}\label{E1}\textbf{x}_{k}-\textbf{x}_{0}\in \textit{K}_{k}(\textit{A}, \textbf{r}_{0})
=span(\textbf{r}_{0}, \textit{A}\textbf{r}_{0},\dots,\textit{A}^{k-1}\textbf{r}_{0}),\end{equation}
\begin{equation}\label{E2}\textbf{r}_{k}=(\textbf{b}-\textit{A}\textbf{x}_{k})\bot\textit{K}_{k}(\textit{A}^{T}, \textbf{y})
=span(\textbf{y},
\textit{A}^T\textbf{y},\dots,\textit{A}^{T^{k-1}}\textbf{y}),\end{equation}
\noindent where $\textit{A}^T$ denotes the transpose of
$\textit{A}$.

\noindent Equation $(\ref{E1})$ leads to,
\begin{equation}
\textbf{x}_{k}-\textbf{x}_{0}=-\alpha_{1}\textbf{r}_{0}-\alpha_{2}\textit{A}\textbf{r}_{0}
- \dots -\alpha_{k}\textit{A}^{k-1}\textbf{r}_{0}.
\end{equation}
\noindent Multiplying both sides by $\textit{A}$ and adding and subtracting $\textbf{b}$ on the left hand side gives
\begin{equation}\label{E5}
\textbf{r}_{k}=\textbf{r}_{0}+\alpha_{1}\textit{A}\textbf{r}_{0}+\alpha_{2}\textit{A}^2\textbf{r}_{0}+
\dots +\alpha_{k}\textit{A}^{k}\textbf{r}_{0}.
\end{equation}
\noindent From $(\ref{E2})$, the orthogonality condition gives

$(\textit{A}^{T^{i}}\textbf{y},\textbf{r}_{k})=0$, for $i=0,\dots,k-1,$

\noindent and, by $(\ref{E5})$, we obtain the following system of linear equations
\begin{eqnarray}\label{E6}\begin{cases}
\alpha_{1}(\textbf{y}, \textit{A}\textbf{r}_{0})+ \dots +
\alpha_{k}(\textbf{y}, \textit{A}^k\textbf{r}_{0})=
-(\textbf{y}, \textbf{r}_{0}),\\
\vdots\\
\alpha_{1}(\textit{A}^{T^{k-1}}\textbf{y}, \textit{A}\textbf{r}_{0})+ \dots +
\alpha_{k}(\textit{A}^{T^{k-1}}\textbf{y}, \textit{A}^k\textbf{r}_{0})=
-(\textit{A}^{T^{k-1}}\textbf{y}, \textbf{r}_{0}).
\end{cases}\end{eqnarray}
\noindent If the determinant of the above system is different from zero then its solution exists and allows to obtain $\textbf{x}_{k}$ and $\textbf{r}_{k}$. Obviously, in practice, solving the above system directly for the increasing value of $k$ is not feasible. We shall now see how to solve this system for increasing values of $k$ recursively.

\noindent If we set
\begin{equation}\label{E7}P_{k}(x)=1+\alpha_{1}x+\dots+\alpha_{k}x^k,\end{equation}
\noindent then we can write from $(\ref{E5})$
\begin{equation}\label{E8}
\textbf{r}_{k}=P_{k}(\textit{A})\textbf{r}_{0}.
\end{equation}
\noindent The polynomials $P_{k}$ are commonly known as the residual polynomials, \cite{93:Brezinski}. Another interpretation of the $P_{k}$ can be found in \cite{87:Cybenko}. Moreover if we set
$c_{i}=(\textit{A}^{T^i}\textbf{y}, \textbf{r}_{0})=(\textbf{y}, \textit{A}^i\textbf{r}_{0})$, $i=0,1,\dots$, and if we define the linear functional $c$ on the space of polynomials by
\begin{equation}
c(x^i)=c_{i},  \mbox{ }  i=0,1,\dots,
\end{equation}
\noindent $c$ is completely determined by the sequence $\{c_k\}$ and $c_k$ is said to be the moment of order $k$, \cite{80:Brezinski}. Now, the system (\ref{E6}) can be written as
\begin{equation}
c(x^iP_{k}(x))=0  \mbox{ for } i=0,\dots,k-1.
\end{equation}
These conditions show that $P_{k}$ is the polynomial of degree at most $k$, normalized by the condition $P_{k}(0)=1$, belonging to a family of FOP's with respect to the linear functional $c$, \cite{92:Brezinski, 80:Brezinski}.

Since the constant term of $P_{k}$ in (\ref{E7}) is $1$, it can be written as
\[P_{k}(x)=1+xR_{k-1}(x)\]
\noindent where $R_{k-1}= \alpha_1+ \alpha_2x+...+\alpha_kx^{k-1}$. Replacing $x$ by $\textit{A}$ in the expression of $P_k$ and multiplying both sides by $\textbf{r}_{0}$ and using $(\ref{E8})$, we get
\[\textbf{r}_{k}=\textbf{r}_{0}+\textit{A}R_{k-1}(A)\textbf{r}_{0},\]
\noindent which can be written as
\[b-\textit{A}\textbf{x}_{k}=b-\textit{A}\textbf{x}_{0}+\textit{A}R_{k-1}(A)\textbf{r}_{0},\]
\[-\textit{A}\textbf{x}_{k}=-\textit{A}\textbf{x}_{0}+\textit{A}R_{k-1}(A)\textbf{r}_{0},\]
\noindent multiplying both sides by $-\textit{A}^{-1}$, we get
\[\textbf{x}_{k}=\textbf{x}_{0}-R_{k-1}(A)\textbf{r}_{0},\]
\noindent which shows that $\textbf{x}_{k}$ can be computed from $\textbf{r}_{k}$ without inverting $\textit{A}$.

\section{Formal orthogonal polynomials}
The orthogonal polynomials $P_{k}$ defined in the previous section are given by the determinantal formula, \cite{97:Brezinski,93:Brezinski}
\begin{equation}
P_{k}(x)=\frac{\left\vert\begin{array}{cccc}
1 & \cdots & x^k\\
c_{0} & \cdots & c_{k}\\
\vdots &  & \vdots\\
c_{k-1} & \cdots & c_{2k-1}\\
\end{array}\right\vert}{\left\vert\begin{array}{cccc}
c_{1} & \cdots & c_{k}\\
\vdots &  & \vdots\\
c_{k} & \cdots & c_{2k-1}\\
\end{array}\right\vert},
\end{equation}
\noindent where the denominator of this polynomial is $\textit{H}^{(1)}_{k}$, \cite{93:Brezinski}. Obviously, $P_{k}$ exists if and only if the Hankel determinant $\textit{H}^{(1)}_{k}\neq0.$ Thus, $P_{k+1}$ exists if and only if $\textit{H}^{(1)}_{k+1}\neq0$.  We assume that $\forall k$, $\textit{H}^{(1)}_{k}\neq0$. If for some $k$, $\textit{H}^{(1)}_{k}=0$, then $P_{k}$ does not exist and breakdown occurs in the algorithm (in practice the breakdown can occur even if $\textit{H}^{(1)}_{k}\approx0$).

Let us now define a linear functional $c^{(1)},$ \cite{93:Brezinski,95:Baheux}, on the space of real polynomials as $c^{(1)}(x^{i})=c(x^{i+1})=c_{i+1}$ and let $P^{(1)}_{k}$ be a family of orthogonal polynomials with respect to $c^{(1)}$. These polynomials are called monic polynomials, \cite{93:Brezinski,95:Baheux}, because their highest degree coefficients are always 1, and are given by the following formula
\begin{equation}
P^{(1)}_{k}(x)=\frac{\left\vert\begin{array}{cccc}
c_{1} & \cdots & c_{k+1}\\
\vdots &    & \vdots\\
c_{k} & \cdots & c_{2k}\\
1 & \cdots & x^k\\
\end{array}\right\vert}{\left\vert\begin{array}{cccc}
c_{1} & \cdots & c_{k}\\
\vdots &  & \vdots\\
c_{k} & \cdots & c_{2k-1}\\
\end{array}\right\vert}.
\end{equation}
\noindent $P^{(1)}_{k}(x)$ also exists if and only if the Hankel determinant $H^{(1)}_{k}\neq0$, \cite{93:Brezinski,95:Baheux}, which is also a condition for the existence of $P_{k}(x)$. There exist many recurrence relations between the two adjacent families of polynomials $P_{k}$ and $P^{(1)}_{k}$, \cite{95:Baheux,93:Brezinski,97:Brezinski,91:Brezinski}. Some of these relations have been reviewed in \cite{94:Zaglia} and studied in details in \cite{94:Baheux,97:Brezinski}. More of these relations have been studied in \cite{10:Farooq}, leading to new Lanczos-type algorithms.

A Lanczos-type algorithm consists in computing $P_{k}$ recursively, then $\textbf{r}_{k}$ and finally $\textbf{x}_{k}$ such that $\textbf{r}_{k}=\textbf{b}-\textit{A}\textbf{x}_{k}$, without inverting $A$. In exact arithmetic, this should give the solution to the system $\textit{A}\textbf{x}=\textbf{b}$ in at most $n$ steps \cite{52:Lanczos,91:Brezinski}, where $n$ is the dimension of the system. For more details, see \cite{93:Brezinski,99:Brezinski}.

\section{Recalling some existing algorithms}
In the following we will recall some of the most recent and efficient Lanczos-type algorithms to be used in the switching framework. The reader should consult the relevant literature for more details.

\subsection{Algorithm $A_{12}$}
Algorithm $A_{12}$ is based on relation $A_{12}$, \cite{10:Farooq}. For details on the derivation of the polynomial $A_{12}$, its coefficients and the algorithm itself, please refer to \cite{10:Farooq}.
\noindent The pseudo-code of Algorithm $A_{12}$ can be described as follows.
\begin{algorithm}[H]
\caption{Algorithm $A_{12}$}
\begin{algorithmic}[1]
\STATE Choose $x_{0}$ and $y$ such that $y\neq0$,
\STATE Choose $\epsilon$ small and positive, as a tolerance,
\STATE Set $r_{0}=b-Ax_{0}$, $y_{0}=y$, $p=Ar_{0}$, $p_{1}=Ap$, $c_{0}=(y, r_{0})$,
\STATE  $c_{1}=(y, p)$, $c_{2}=(y, p_{1})$, $c_{3}=(y, Ap_{1})$, $\delta=c_{1}c_{3}-c_{2}^2$,
\STATE $\alpha=\frac{c_{0}c_{3}-c_{1}c_{2}}{\delta}$, $\beta=\frac{c_{0}c_{2}-c_{1}^2}{\delta}$,
$r_{1}=r_{0}-\frac{c_{0}}{c_{1}}p$, $x_{1}=x_{0}+\frac{c_{0}}{c_{1}}r_{0}$,
\STATE $r_{2}=r_{0}-\alpha p+\beta p_{1}$, $x_{2}=x_{0}+\alpha r_{0}-\beta p$,
\STATE $y_{1}=A^{T}y_{0}$, $y_{2}=A^{T}y_{1}$, $y_{3}=A^{T}y_{2}$.
\FOR {k = 3, 4,\dots, n} \label{k=3}
\STATE $y_{k+1}=A^{T}y_{k}$, $q_{1}=Ar_{k-1}$, $q_{2}=Aq_{1}$, $q_{3}=Ar_{k-2}$,
\STATE$a_{11}=(y_{k-2}, r_{k-2})$, $a_{13}=(y_{k-3}, r_{k-3})$, $a_{21}=(y_{k-1}, r_{k-2})$, $a_{22}=a_{11}$,
\STATE $a_{23}=(y_{k-2}, r_{k-3})$, $a_{31}=(y_{k}, r_{k-2})$,$a_{32}=a_{21}$, $a_{33}=(y_{k-1}, r_{k-3})$,
\STATE $s=(y_{k+1}, r_{k-2})$, $t=(y_{k}, r_{k-3})$,$F_{k}=-\frac{a_{11}}{a_{13}}$,
\STATE $b_{1}=-a_{21}-a_{23}F_{k}$, $b_{2}=-a_{31}-a_{33}F_{k}$, $b_{3}=-s-tF_{k}$,
\STATE $\Delta_{k}=a_{11}(a_{22}a_{33}-a_{32}a_{23})+a_{13}(a_{21}a_{32}-a_{31}a_{22})$,
\STATE $B_{k}=\frac{b_{1}(a_{22}a_{33}-a_{32}a_{23})+a_{13}(b_{2}a_{32}-b_{3}a_{22})}{\Delta_{k}}$,
\STATE $G_{k}=\frac{b_{1}-a_{11}B_{k}}{a_{13}}$, $C_{k}=\frac{b_{2}-a_{21}B_{k}-a_{23}G_{k}}{a_{22}}$, $A_{k}=\frac{1}{C_{k}+G_{k}}$,
\STATE $r_{k}=A_{k}\{q_{2}+B_{k}q_{1}+C_{k}r_{k-2}+F_{k}q_{3}+G_{k}r_{k-3}\}$,
\STATE $x_{k}=A_{k}\{C_{k}x_{k-2}+G_{k}x_{k-3}-(q_{1}+B_{k}r_{k-2}+F_{k}r_{k-3})\}$,
\IF {$||r_{k}|| \leq \epsilon $,}
\STATE $x = x_{k}$, Stop.
\ENDIF
\ENDFOR
\end{algorithmic}
\end{algorithm}

\subsection{Algorithm $A_{4}$}Algorithm $A_{4}$ is based on relation $A_{4}$. Its pseudo-code is as follows. For more details see \cite{94:Baheux,95:Baheux}.
\begin{algorithm}[H]
\caption{Algorithm $A_{4}$}
\begin{algorithmic}[1]
\STATE Choose $x_{0}$ and $y$ such that $y\neq0$,
\STATE Choose $\epsilon$ small and positive as a tolerance,
\STATE Set $r_{0}=b-Ax_{0}$, $y_{0} = y$,
\FOR {k = 0, 1,\dots, n}\label{k=0}
\STATE $E_{k+1}=-\frac{(y_{k}, r_{k})}{(y_{k-1}, r_{k-1})}$, for $k\geq1$, and $E_{1}=0$,
\STATE $B_{k+1}=-\frac{(y_{k}, Ar_{k})-E_{k+1}(y_{k}, r_{k-1})}{(y_{k}, r_{k})}$,
\STATE $A_{k+1}=\frac{1}{B_{k+1}+E_{k+1}}$,
\STATE $x_{k+1}=A_{k+1}\{B_{k+1}x_{k}+E_{k+1}x_{k-1}-r_{k}\}$,
\STATE $r_{k+1}=A_{k+1}\{Ar_{k}+B_{k+1}r_{k}+E_{k+1}r_{k-1}\}$.
\IF {$||r_{k+1}|| \leq \epsilon $,}
\STATE $y_{k+1}=A^Ty_{k}$,
\ENDIF\ENDFOR
\end{algorithmic}
\end{algorithm}

\subsection{Algorithm $A_{5}/B_{10}$}
Algorithm $A_{5}/B_{10}$ is based on relations $A_{5}$ and $B_{10}$, first investigated in \cite{94:Baheux,95:Baheux}. Its pseudo-code is as follows.
\begin{algorithm}[H]
\caption{Algorithm $A_{5}/B_{10}$}
\begin{algorithmic}[1]
\STATE Choose $x_{0}$, $y$ and tolerance $\epsilon \geq 0 $;
\STATE Set $r_{0}=b-Ax_{0}$, $p_{0} = r_{0}$, $y_{0} = y$,
\STATE $A_{1} = -\frac{(y_{0}, r_{0})}{(y_{0}, Ar_{0})}$, $C^1_{0}= 1$,
\STATE $r_{1} = r_{0}+A_{1}Ar_{0}$, $x_{1} = x_{0}-A_{1}r_{0}$.
\FOR {k = 1,2,3,\dots,n}
\STATE $y_{k} = A^Ty_{k-1}$,
\STATE $D_{k+1} = -\frac{(y_{k}, r_{k})}{C^1_{k-1}(y_{k}, p_{k-1})}$,
\STATE $p_{k}=r_{k}+D_{k-1}C^1_{k-1}p_{k-1}$
\STATE $A_{k+1} = -\frac{(y_{k}, r_{k})}{(y_{k}, Ap_{k})}$,
\STATE $r_{k+1} = r_{k}+A_{k+1}Ap_{k}$,
\STATE $x_{k+1} = x_{k}-A_{k+1}p_{k}$.
\IF {$||r_{k+1}|| \neq \epsilon$ and $A_{k} \neq \epsilon$,}
\STATE $C^1_{k} = \frac{C^1_{k-1}}{A_{k}}$.
\ENDIF
\ENDFOR
\end{algorithmic}
\end{algorithm}

\subsection{Algorithm $A_{8}/B_{10}$}
The pseudo-code of $A_{8}/B_{10}$, \cite{94:Baheux,95:Baheux}, is as follows.
\begin{algorithm}[H]
\caption{Algorithm $A_{8}/B_{10}$}
\begin{algorithmic}[1]
\STATE Choose $x_{0}$ and $y$ such that $y\neq0$.\\
\STATE Set $r_0 = b - Ax_0$,\\
\STATE $z_{0}=r_{0}$,\\
\STATE $y_{0}=y$,\\
\FOR {$k=0,1,2,\dots$,n}
\STATE $A_{k+1}=-\frac{(\textbf{y}_{k}, \textbf{r}_{k})}{(\textbf{y}_{k}, \textit{A}\textbf{z}_{k})}$,\\
\STATE $r_{k+1}=r_{k}+A_{k+1}Az_{k}$,\\
\STATE $x_{k+1}=x_{k}-A_{k+1}z_{k}$.\\
\IF{$||r_{k+1}|| \neq \epsilon$,}
\STATE $y_{k+1}=A^{T}y_{k}$,\\
\STATE $C^1_{k+1}=\frac{1}{A_{k+1}}$,\\
\STATE $B^{1}_{k+1}=-\frac{C^{1}_{k+1}(y_{k+1}, r_{k+1})}{(y_{k}, Az_{k})}$,\\
\STATE $z_{k+1}=B^{1}_{k+1}z_{k}+C^{1}_{k+1}r_{k+1}$.
\ENDIF
\ENDFOR
\end{algorithmic}
\end{algorithm}

\section{Switching between algorithms to avoid breakdown}
When a Lanczos-type algorithm fails, this is due to the non-existence of some coefficients of the recurrence relations on which the algorithm is based. The iterate which causes these coefficients not to exist does not cause and should not necessarily cause any problems when used in another Lanczos-type algorithm, based on different recurrence relations. It is therefore obvious that one may consider switching to this other algorithm, when breakdown occurs. This allows the algorithm to work in a Krylov space with a different basis. It is therefore also possible to remedy breakdown by switching. Note that restarting the same algorithm after a pre-set number of iterations works well too, \cite{12:Farooq}

\newpage
\subsection{\bf{Switching strategies}}
Different strategies can be adopted for switching between two or more algorithms. These are as follows.

\begin{enumerate}
\item \textbf{ST1: Switching after breakdown:}
Start a particular Lanczos algorithm until a breakdown occurs, then switch to another Lanczos algorithm, initializing the latter with the last iterate of the failed algorithm. We call this strategy ST1.

\item \textbf{ST2: Pre-emptive switching:}
Run a Lanczos-type algorithm for a fixed number of iterations, halt it and then switch to another Lanczos-type algorithm, initializing it with the last iterate of the first algorithm. Note that there is no way to guarantee that breakdown would not occur before the end of the interval. This strategy is called ST2.

\item \textbf{ST3: Breakdown monitoring:} Provided monotonicity of reduction in the absolute value of the denominators in the coefficients of the polynomials involved can be established, breakdown can be monitored as follows. Evaluate regularly those coefficients with denominators that are likely to become zero. Switch to another algorithm when the absolute value of any of these denominators drops below a certain level. This is strategy ST3.
\end{enumerate}

\subsection{\bf{A generic switching algorithm}}
Suppose we have a set of Lanczos-type algorithms and we want to switch from one algorithm to another using one of the above mentioned strategies ST1, ST2 or ST3.
\begin{algorithm}[H]
\caption{Generic switching algorithm}
\begin{algorithmic}[1]
\STATE Start the most stable algorithm, if known.
\STATE Choose a switching strategy from $\{$\bf{ST1, ST2, ST3}$\}$.
\IF{\textbf{ST1}}
\STATE Continue with current algorithm until it halts;
\IF{solution is obtained}
\STATE Stop.
\ELSE
\STATE switch to another algorithm;
\STATE initialize it with current iterate;
\STATE Go to 4.
\ENDIF
\ELSIF{\textbf{ST2}}
\STATE Continue with current algorithm for a fixed number of iterations until it stops;
\IF{solution is obtained}
\STATE Stop.
\ELSE
\STATE switch to another algorithm,
\STATE initialize it with the current iterate,
\STATE Go to 13.
\ENDIF
\ELSE
\STATE Continue with current algorithm and monitor certain parameters for breakdown, until it halts;
\IF{solution is obtained} 
\STATE Stop.
\ELSE
\STATE switch to another algorithm,
\STATE initialize it with the current iterate,
\STATE Go to 22.
\ENDIF
\ENDIF
\end{algorithmic}
\end{algorithm}

\noindent However, it is important to mention that we have considered only ST2 in this paper. The convergence tolerance in all of the tests performed is $\epsilon = 1.0e^{-013}$ and the number of iterations per cycle is fixed to $20$.

\subsubsection{\bf{Switching between algorithms $A_{4}$ and $A_{12}$}}
In the following, we start with either $A_4$ or $A_{12}$, run it for a fixed number of iterations (cycle) chosen arbitrarily, before switching to the other. The results of this switching algorithm, are compared to those obtained with algorithms $A_{4}$ and $A_{12}$ run individually. We are not changing any of the parameters involved in both algorithms. Details of $A_4$ can be found in \cite{94:Baheux}.
\begin{algorithm}[H]
\caption{Switching between $A_{4}$ and $A_{12}$}
\setstretch{0.5}
\begin{algorithmic}[1]
\STATE Choose $x_{0}$ and $y$ such that $y\neq0$,
\STATE set $r_{0}=b-Ax_{0}$, $y_{0}=y$,
\STATE start either algorithm,
\STATE run current algorithm for a fixed number of iterations (a cycle) or until it halts;
\IF{solution is obtained}
\STATE stop;
\ELSE
\STATE switch to the algorithm not yet run;
\STATE initialize it with the current iterate;
\STATE go to 4;
\ENDIF
\end{algorithmic}
\end{algorithm}
\noindent {\bf Remark:} Since restarting can be just as effective as switching, it is easier to implement a random choice between $A_{4}$ and $A_{12}$ at the end of every cycle. Let heads be $A_{4}$ and tails be $A_{12}$. At the toss of a coin, if it shows heads and the algorithm running in the last cycle was $A_{4}$, then the switch is a restart. If the coin shows tails then the switch is a ``proper" switch, and $A_{12}$ is called upon. In the numerical results presented below, this is what has been implemented. For more details about restarting see, \cite{12:Farooq}.

\subsubsection{\bf{Switching between $A_{4}$ and $A_{5}/B_{10}$ algorithm}}
Start with $A_{5}/B_{10}$, (details of $A_5/B_{10}$ can be found in \cite{94:Baheux,95:Baheux}) do a few iterations and then switch to either $A_{4}$ or $A_{5}/B_{10}$. The procedure is as Algorithm $4$ below.
\begin{algorithm}[H]
\caption{Switching between $A_{4}$ and $A_{5}/B_{10}$}
\begin{algorithmic}[1]
\STATE Choose $x_{0}$ and $y$ such that $y \neq 0$;
\STATE set $r_{0} = b-Ax_{0}$, $y_{0} = y$, $p_{0} = r_{0}$;
\STATE start with either $A_4$ or $A_{5}/B_{10}$;
\STATE run it for a fixed number of iterations (cycles) or until it halts
\IF{solution is obtained}
\STATE stop;
\ELSE
\STATE switch to either $A_{4}$ or $A_{5}/B_{10}$; initialize it with the last iterate of the algorithm running in the last cycle;
\STATE go to 4;
\ENDIF
\end{algorithmic}
\end{algorithm}

\subsubsection{\bf{Switching between $A_{4}$ and $A_{8}/B_{10}$}}
Start with either $A_{8}/B_{10}$ (details of $A_8/B_{10}$ can be found in \cite{94:Baheux,95:Baheux}) or $A_4$; do a few iterations and then switch to either of them chosen randomly. If the chosen algorithm happens to be the same as the one running in the last cycle, then it is a case of restarting. Otherwise, it is switching. The algorithm is as follows.
\begin{algorithm}[H]
\caption{Switching between $A_{4}$ and $A_{8}/B_{10}$}
\setstretch{0.5}
\begin{algorithmic}[1]
\STATE Choose $x_{0}$ and $y$ such that $y \neq 0$;
\STATE set $r_{0}=b-Ax_{0}$, $y_{0}=y$, $p_{0}=r_{0}$;
\STATE start either $A_{4}$ or $A_{8}/B_{10}$;
\STATE run it for a fixed number of iterations (cycle), or until it halts;
\IF{solution is obtained} 
\STATE stop;
\ELSE
\STATE switch to either $A_{4}$ or $A_{8}/B_{10}$;
\STATE initialize it with the iterate of the algorithm run in the last cycle;
\STATE go to 4.
\ENDIF
\end{algorithmic}
\end{algorithm}

\subsubsection{\bf{Switching between $A_{5}/B_{10}$ and $A_{8}/B_{10}$}}
Here again, switching and restarting are combined in a random way. Start with either $A_{8}/B_{10}$ or $A_{5}/B_{10}$. After a pre-set number of iterations (cycle), switch to either $A_{5}/B_{10}$ or $A_{8}/B_{10}$, randomly chosen. If the chosen algorithm to switch to is the same as the one running in the last cycle then we a have a case of restarting; else it is switching. The algorithm is as follows.
\begin{algorithm}[H]
\caption{Switching between $A_{5}/B_{10}$ and $A_{8}/B_{10}$}
\setstretch{0.5}
\begin{algorithmic}[1]
\STATE Choose $x_{0}$ and $y$ such that $y\neq0$;
\STATE set $r_{0}=b-Ax_{0}$, $y_{0}=y$, $z_{0}=r_{0}$;
\STATE start either $A_{8}/B_{10}$ or $A_{5}/B_{10}$;
\STATE run it for a fixed number of iterations;
\IF{solution is not found} 
\STATE halt current algorithm;
\STATE switch to either $A_{5}/B_{10}$ or $A_{8}/B_{10}$;
\STATE initialize it with the last iterate of the algorithm running in the last cycle;\\
\STATE go to 4; \\
\ELSE
\STATE solution found; stop;
\ENDIF
\end{algorithmic}
\end{algorithm}

\subsubsection{\bf{Numerical results}}
Algorithms $1$, $2$, $3$, $4$, \cite{94:Baheux,95:Baheux,10:Farooq} and Algorithms 6, 7, 8 and 9, \cite{10:Farooq} have been implemented in Matlab and applied to a number of small to medium size problems. The test problems we have used arise in the 5-point discretisation of the operator $-\frac{\partial^{2}}{\partial x^2}-\frac{\partial^{2}}{\partial y^2}+\gamma\frac{\partial}{\partial x}$ on a rectangular region \cite{94:Baheux,95:Baheux}. Comparative results are obtained on instances of the problem $\textit{A}\textbf{x}=\textbf{b}$ with $\textit{A}$ and $\textbf{b}$ as below, and with dimensions of $A$ and $\textbf{b}$ ranging from $n=10$ to $n=100$.
\[A=\left(\begin{array}{ccccccc}
B & -I & \cdots & \cdots  & 0\\
-I & B & -I  &  & \vdots\\
\vdots & \ddots & \ddots & \ddots & \vdots\\
\vdots &  & -I & B & -I\\
0 & \cdots & \cdots & -I & B\\
\end{array}
\right),\]\noindent with
\[B=\left(\begin{array}{ccccccc}
4 & \alpha & \cdots & \cdots & 0\\
\beta & 4 & \alpha &  & \vdots\\
\vdots & \ddots & \ddots & \ddots & \vdots\\
\vdots & & \beta & 4 & \alpha\\
0 & \cdots &  & \beta & 4\\
\end{array}\right),\]
\noindent and $\alpha=-1+\delta$, $\beta=-1-\delta$.  The parameter $\delta$ takes the values $0.0$, $0.2$, $5$ and $8$ respectively. The right hand side $\textbf{b}$ is taken to be $\textbf{b}=\textit{A}\textbf{X}$, where $\textbf{X}=(1, 1, \dots, 1)^{T}$, is the solution of the system. The dimension of $\textit{B}$ is $10$. When $\delta=0$, the coefficient matrix $A$ is symmetric and the problem is easy to solve because the region is a regular mesh, \cite{06:Gérard}. For all other values of $\delta$, the matrix $A$ is non-symmetric and the problem is comparatively harder to solve as the region is not a regular mesh.

\subsubsection{\bf{Numerical results}}

Results obtained with Algorithms 1, 2, 3, 4, and Algorithms 6, 7, 8 and 9 on Baheux-type problems of different dimensions, for different values of $\delta$ are presented in tables 1, 2, 3 and 4, below. Algorithms $1$, $2$, $3$, $4$, executed individually, could only solve problems of dimensions 40 or below. In contrast, the switching algorithms, Algorithms 6, 7, 8 and 9, solved all problems of dimensions up to 4000. These results show that the switching algorithms are far superior to any one of the algorithms considered individually.
These echo those obtained by restarting the same algorithm after a predefined number of iterations, \cite{12:Farooq}.
\begin{table}[H]
\caption{Numerical results for $\delta=0$} 
\vspace{0.5cm}
\centering 
\scalebox{0.60}{
\begin{tabular}{|c|c|c|c|c|c|c|c|c|c|c|c|c|}
\hline
Dim of Prob & \multicolumn{2}{|c|}{Algorithm 6} &  \multicolumn{2}{|c|}{Algorithm 7} &  \multicolumn{2}{|c|}{Algorithm 8} &  \multicolumn{2}{|c|}{Algorithm 9}\\ [0.5ex] 
\hline\hline
$n$ & $||r_{k}||$  & T(s) & $||r_{k}||$ & T(s) & $||r_{k}||$ & T(s) & $||r_{k}||$ & T(s)\\[0.5ex] 
\hline
20 & $5.5067e^{-014}$ & 0.0012 & $3.8545e^{-014}$ & 0.0010 & $7.4781e^{-014}$ & 0.0040 & $8.0533e^{-014}$ & 0.0018 \\
40 & $7.5417e^{-014}$ & 0.0041 & $4.5076e^{-014}$ & 0.0053 & $8.9208e^{-014}$ & 0.0057 & $7.2481e^{-014}$ & 0.0040\\
60 & $9.6638e^{-014}$ & 0.0057 & $2.5330e^{-014}$ & 0.0057 & $7.5107e^{-014}$ & 0.0085 & $5.8162e^{-014}$ & 0.0067\\
80 & $9.9082e^{-014}$ & 0.0075 & $6.0185e^{-014}$ & 0.0071 & $7.0866e^{-014}$ & 0.0088 & $5.7266e^{-014}$ & 0.0101\\
100 & $2.1487e^{-014}$ & 0.0095 & $2.5839e^{-014}$ & 0.0078 & $8.0262e^{-014}$ & 0.0098 & $8.1373e^{-014}$ & 0.0100\\
200 & $7.4236e^{-014}$ & 0.0723 & $9.9667e^{-014}$ & 0.0159 & $7.9045e^{-014}$ & 0.0337 & $9.1830e^{-014}$ & 0.0352\\
400 & $7.7419e^{-014}$ & 0.0661 & $8.5151e^{-014}$ & 0.2156 & $9.7418e^{-014}$ & 0.2243 & $9.4697e^{-014}$ & 0.2315\\
600 & $9.0290e^{-014}$ & 0.0794 & $7.9373e^{-014}$ & 0.4735 & $9.9269e^{-014}$ & 1.9625 & $9.4307e^{-014}$ & 0.7457\\
800 & $9.2116e^{-014}$ & 0.5660 & $9.5227e^{-014}$ & 0.9395 & $7.7294e^{-014}$ & 3.0326 & $7.7356e^{-014}$ & 1.4319\\
1000 & $8.8463e^{-014}$ & 0.8509 & $8.6238e^{-014}$ & 2.0539 & $9.9181e^{-14}$ & 4.3479 & $9.4512e^{-014}$ & 2.8984\\
2000 & $9.7242e^{-014}$ & 4.8079 & $9.1973e^{-014}$ & 8.6364 & $8.1319e^{-014}$ & 12.8696 & $9.1193e^{-014}$ & 12.8696\\
3000 & $9.6993e^{-014}$ & 9.8130 & $7.7507e^{-014}$ & 16.5795 & $9.7827e^{-014}$ & 22.3386 & $8.2725e^{-014}$ & 22.3386\\
4000 & $9.2641e^{-014}$ & 17.4673 & $8.9681e^{-014}$ & 23.7658 & $9.8438e^{-014}$ & 44.4430 & $9.7911e^{-014}$ & 44.4567\\
\hline
\end{tabular}}
\end{table}

\begin{table}[H]
\caption{Numerical results for $\delta=0.2$} 
\vspace{0.5cm}
\centering 
\scalebox{0.60}{
\begin{tabular}{|c|c|c|c|c|c|c|c|c|c|c|c|c|}
\hline
Dim of Prob & \multicolumn{2}{|c|}{Algorithm 6} &  \multicolumn{2}{|c|}{Algorithm 7} &  \multicolumn{2}{|c|}{Algorithm 8} &  \multicolumn{2}{|c|}{Algorithm 9}\\ [0.5ex] 
\hline\hline
$n$ & $||r_{k}||$  & T(s) & $||r_{k}||$ & T(s) & $||r_{k}||$ & T(s) & $||r_{k}||$ & T(s)\\[0.5ex] 
\hline
20 & $6.0041e^{-014}$ & 0.0044 & $1.5104e^{-014}$ & 0.0029 & $1.8618e^{-014}$ & 0.0056 & $8.0533e^{-014}$ & 0.0044\\
40 & $1.9868e^{-014}$ & 0.0064 & $4.6814e^{-014}$ & 0.0068 & $3.4094e^{-014}$ & 0.0089 & $7.3581e^{-014}$ & 0.0100\\
60 & $6.0788e^{-014}$  & 0.0134 & $6.3548e^{-014}$ & 0.0104 & $2.9827e^{-014}$ & 0.0113 & $9.6583e^{-014}$ & 0.0262\\
80 & $8.8550e^{-014}$  & 0.0159 & $9.5483e^{-014}$ & 0.0108 & $8.4187e^{-014}$ & 0.0120 & $7.0744e^{-014}$ & 0.0269\\
100 & $5.8020e^{-014}$  & 0.0144 & $6.6962e^{-014}$ & 0.0151 & $7.3889e^{-014}$ & 0.0126 & $7.5236e^{-014}$ & 0.0273\\
200 & $9.0970e^{-014}$  & 0.0213 & $9.8054e^{-014}$ & 0.0353 & $8.6331e^{-014}$ & 0.0313 & $8.1282e^{-014}$ & 0.0352\\
400 & $6.5593e^{-014}$  & 0.0748 & $9.3591e^{-014}$ & 0.1054 & $6.6660e^{-014}$ & 0.1875 & $8.4316e^{-014}$ & 0.2315\\
600 & $9.6153e^{-014}$  & 0.1802 & $8.8169e^{-014}$ & 0.6066 & $6.8135e^{-014}$ & 0.6751 & $5.9937e^{-014}$ & 0.7457\\
800 & $9.8605e^{-014}$  & 0.5922 & $8.8399e^{-014}$ & 0.9088 & $8.2550e^{-014}$ & 1.1436 & $7.0295e^{-014}$ & 1.4319\\
1000 & $9.7823e^{-014}$  & 0.8222 & $7.7898e^{-014}$ & 1.3020 & $7.4540e^{-014}$ & 2.1302 & $7.7204e^{-014}$ & 2.8984\\
2000 & $7.9753e^{-014}$  & 4.3416 & $9.4241e^{-014}$ & 4.7668 & $9.5282e^{-014}$ & 10.5976 & $9.1570e^{-014}$ & 8.5787\\
3000 & $8.7448e^{-014}$  & 10.0287 & $9.6831e^{-014}$ & 12.1561 & $9.9608e^{-014}$ & 25.1173 & $9.2806e^{-014}$ & 21.7380\\
4000 & $5.8412e^{-014}$  & 12.9390 & $9.7580e^{-014}$ & 23.3993 & $9.9270e^{-014}$ & 38.9051 & $9.7911e^{-014}$ & 39.0164\\
\hline
\end{tabular}}
\end{table}

\begin{table}[H]
\caption{Numerical results for $\delta=5$} 
\vspace{0.5cm}
\centering 
\scalebox{0.60}{
\begin{tabular}{|c|c|c|c|c|c|c|c|c|c|c|c|c|}
\hline
Dim of Prob & \multicolumn{2}{|c|}{Algorithm 6} &  \multicolumn{2}{|c|}{Algorithm 7} &  \multicolumn{2}{|c|}{Algorithm 8} &  \multicolumn{2}{|c|}{Algorithm 9}\\ [0.5ex] 
\hline\hline
$n$ & $||r_{k}||$  & T(s) & $||r_{k}||$ & T(s) & $||r_{k}||$ & T(s) & $||r_{k}||$ & T(s)\\[0.5ex] 
\hline
20 & $2.5092e^{-014}$ & 0.0079 & $9.1438e^{-014}$ & 0.0060 & $3.5839e^{-014}$ & 0.0067 & $8.0533e^{-014}$ & 0.0038\\
40 & $7.4721e^{-014}$ & 0.0202 & $1.9575e^{-014}$ & 0.0171 & $8.9026e^{-014}$ & 0.0079 & $2.9089e^{-014}$ & 0.0078\\
60 & $9.1477e^{-014}$  & 0.0232 & $5.4122e^{-014}$ & 0.0218 & $3.2486e^{-014}$ & 0.0207 & $5.6811e^{-014}$ & 0.0103\\
80 & $2.5165e^{-014}$  & 0.0275 & $7.9597e^{-014}$ & 0.0277 & $7.7553e^{-014}$ & 0.0255 & $5.7266e^{-014}$ & 0.0101\\
100 & $8.9244e^{-014}$  & 0.0315 & $8.4269e^{-014}$ & 0.0295 & $3.2898e^{-014}$ & 0.0308 & $8.1373e^{-014}$ & 0.0100\\
200 & $8.5274e^{-014}$  & 0.0410 & $8.3743e^{-014}$ & 0.0402 & $4.5501e^{-014}$ & 0.0499 & $9.1830e^{-014}$ & 0.0352\\
400 & $9.5005e^{-014}$  & 0.0965 & $3.3013e^{-014}$ & 0.2662 & $4.3032e^{-014}$ & 0.1856 & $9.4697e^{-014}$ & 0.2315\\
600 & $9.3474e^{-014}$  & 0.2318 & $2.7456e^{-014}$ & 0.7717 & $8.1621e^{-014}$ & 0.6303 & $9.4307e^{-014}$ & 0.7457\\
800 & $7.2197e^{-014}$  & 0.6875 & $9.3718e^{-014}$ & 0.8720 & $9.2023e^{-014}$ & 1.0426 & $7.7356e^{-014}$ & 1.4319\\
1000 & $9.4690e^{-014}$  & 1.7006 & $8.2225e^{-014}$ & 2.4118 & $6.7618e^{-014}$ & 2.6779 & $9.4512e^{-014}$ & 4.2678\\
2000 & $7.0752e^{-014}$  & 9.2566 & $8.8127e^{-014}$ & 6.9938 & $4.6266e^{-014}$ & 11.2416 & $9.1193e^{-014}$ & 11.0604\\
3000 & $8.0276e^{-014}$  & 15.5897 & $8.8194e^{-014}$ & 18.8125 & $4.3762e^{-014}$ & 25.3675 & $8.2725e^{-014}$ & 24.6007\\
4000 & $ 9.7667e^{-014}$  & 29.7400 & $8.9260e^{-014}$ & 30.4619 & $8.5908e^{-014}$ & 42.5710 & $9.7911e^{-014}$ & 41.5276\\
\hline
\end{tabular}}
\end{table}

\begin{table}[H]
\caption{Numerical results for $\delta=8$} 
\vspace{0.5cm}
\centering 
\scalebox{0.60}{
\begin{tabular}{|c|c|c|c|c|c|c|c|c|c|c|c|c|}
\hline
Dim of Prob & \multicolumn{2}{|c|}{Algorithm 6} &  \multicolumn{2}{|c|}{Algorithm 7} &  \multicolumn{2}{|c|}{Algorithm 8} &  \multicolumn{2}{|c|}{Algorithm 9}\\ [0.5ex] 
\hline\hline
$n$ & $||r_{k}||$  & T(s) & $||r_{k}||$ & T(s) & $||r_{k}||$ & T(s) & $||r_{k}||$ & T(s)\\[0.5ex] 
\hline
20 & $8.1056e^{-014}$ & 0.0069 & $8.7622e^{-017}$ & 0.0049 & $3.1380e^{-014}$ & 0.0075 & $8.0533e^{-014}$ & 0.0048\\
40 & $9.1880e^{-014}$ & 0.0243 & $8.1258e^{-014}$ & 0.0061 & $8.2974e^{-014}$ & 0.0216 & $2.9089e^{-014}$ & 0.0078\\
60 & $7.0558e^{-014}$ & 0.0299 & $8.3090e^{-014}$ & 0.0220 & $9.9891e^{-014}$ & 0.0284 & $5.6811e^{-014}$  & 0.0103\\
80 & $9.8855e^{-014}$ & 0.0346 & $8.5600e^{-014}$ & 0.0287 & $9.7330e^{-014}$ & 0.0282 & $5.7266e^{-014}$  & 0.0120\\
100 & $8.7391e^{-014}$ & 0.0382 & $8.7752e^{-014}$ & 0.0303 & $8.5960e^{-014}$ & 0.0363 & $8.1373e^{-014}$ & 0.0137\\
200 & $9.5793e^{-014}$ & 0.0700 & $4.3407e^{-014}$ & 0.0457 & $9.8381e^{-014}$ & 0.0708 & $9.1830e^{-014}$ & 0.0352\\
400 & $6.0799e^{-014}$ & 0.1292 & $9.6421e^{-014}$ & 0.2399 & $9.6706e^{-014}$ & 0.3125 & $9.4697e^{-014}$ & 0.2315\\
600 & $9.6186e^{-014}$ & 0.3432 & $8.5386e^{-014}$ & 0.5535 & $8.9805e^{-014}$ & 0.9279 & $9.4307e^{-014}$ & 0.7457\\
800 & $8.8932e^{-014}$ & 0.6942 & $3.1458e^{-014}$ & 1.3329 & $7.5301e^{-014}$ & 1.0612 & $7.7356e^{-014}$ & 1.4319\\
1000 & $9.7821e^{-014}$ & 1.7060 & $4.9703e^{-014}$ & 2.7150 & $9.6384e^{-014}$ & 2.1978 & $9.4512e^{-014}$ & 2.8984\\
2000 & $7.3843e^{-014}$ & 11.2436 & $8.1578e^{-014}$ & 13.0654 & $8.2557e^{-014}$ & 10.7977 & $9.1193e^{-014}$ & 13.2915\\
3000 & $9.5905e^{-014}$ & 20.0131 & $7.1928e^{-014}$ & 20.9822 & $5.3725e^{-014}$ & 25.3714 & $8.2725e^{-014}$ & 28.4232\\
4000 & $2.1552e^{-014}$ & 31.2356 & $9.4300e^{-014}$ & 40.2119 & $9.3869e^{-014}$ & 40.1782 & $9.7911e^{-014}$ & 45.1523\\
\hline
\end{tabular}}
\end{table}


\section{Conclusion}

We have implemented $A_4$, $A_{12}$, $A_5/B_{10}$ and $A_8/B_{10}$ to solve a number of problems of the type described in Section $5.2.5$ with dimensions ranging from 20 to 4000. The results are compared against those obtained by the switching algorithms, Algorithms 6, 7, 8 and 9 on the same problems. These results show that $A_4$, $A_{12}$, $A_5/B_{10}$ and $A_8/B_{10}$ are not as robust as the switching algorithms. In fact, individual algorithms solved only problems of dimension $n\leq 40$ and that with a poor accuracy. The switching algorithms, however, solved them all with a higher precision. The numerical evidence is strongly in favour of switching.

Based on the above results, it is clear that switching is an effective way to deal with the breakdown in Lanczos-type algorithms. It is also clear that the switching algorithms are more efficient particularly for large dimension problems.

The cost of switching, in terms of CPU time, in ST2 at least, is not substantial, compared to that of the individual algorithms. It is also quite easy to see that it would not be substantial in ST1 since the cost would be similar to that of ST2. Even in the case of monitoring the coefficients that can vanish, i.e. ST3, the cost should only be that of a test of the form: \\ {\bf if}  $|$denominator value$|$ $\leq$ tolerance {\bf then} stop.\\ We have not measured its impact on the overall computing time, but it should not be excessive. This means that switching strategies are worthwhile considering to enhance the efficiency of Lanczos-type algorithms and not just their robustness.

Having said that, further research and experimentation are necessary, particularly on the very large scale instances of SLE's, to establish the superiority of switching algorithms against the state-of-the-art Lanczos-type algorithms with in-built precautions to avoid breakdown such as MRZ and BSMRZ, \cite{94:Baheux,91:Brezinski,na1,99:Brezinski}. Note that these algorithms are attractive for other reasons too, namely their simplicity and easy implementation. This is the subject of on-going research work.

\bibliography{lanczos1}

\noindent Muhammad Farooq,\\
Department of Mathematics, \\
University of Peshawar, \\
Khyber Pakhtunkhwa, 25120, Pakistan \\
mfarooq@upesh.edu.pk, Tel: 00 92 91 9221038, Fax: 00 92 91 9216470 \\

\bigskip

\noindent Abdellah Salhi,\\
Department of Mathematical Sciences, \\
The University of Essex, Wivenhoe Park, \\
Colchester CO4 3SQ, UK \\
as@essex.ac.uk, Tel: 00 44 1206 873022, Fax: 00 44 1206 873043 \\

\end{document}